\def\Statusstring{
            }
            
\documentclass[12pt,letter]{article}

\usepackage{mystyle}
\usepackage{mythmstyle_report}

\begin{document}
\title{Edge-transitivity of Cayley graphs generated by transpositions}

\author{Ashwin~Ganesan%
  \thanks{53 Deonar House, Deonar Village Road, Deonar, Mumbai 400088, Maharashtra, India. Email:  
\texttt{ashwin.ganesan@gmail.com}.}
}

\date{}
\vspace{10cm}

\maketitle

\vspace{-5.5cm}
\begin{flushright}
  \texttt{\Statusstring}\\[1cm]
\end{flushright}
\vspace{+3.0cm}

\begin{abstract}
\noindent  Let $S$ be a set of transpositions generating the symmetric group $S_n$.  The transposition graph of $S$ is defined to be the graph with vertex set $\{1,\ldots,n\}$, and with vertices $i$ and $j$ being adjacent in $T(S)$ whenever $(i,j) \in S$.  In the present note, it is proved that two transposition graphs are isomorphic if and only if the corresponding two Cayley graphs are isomorphic.  It is also proved that the transposition graph $T(S)$ is edge-transitive if and only if the Cayley graph $\Cay(S_n,S)$ is edge-transitive.
\end{abstract}


\vskip 0.3in

\tableofcontents

\vspace{+0.5cm}
\section{Introduction}

Let $X=(V,E)$ be a simple, undirected graph.  An automorphism of $X$ is a permutation of the vertex set that preserves adjacency.  The automorphism group of $X$, denoted by $\Aut(X)$, is the set of all automorphisms of the graph $X$, that is, $\Aut(X) := \{g \in \Sym(V): E^g = E\}$.  A graph $X$ is said to be vertex-transitive if for any two vertices $u, v \in V(X)$, there exists an automorphism $g \in \Aut(X)$ that takes $u$ to $v$.  A graph $X$ is said to be edge-transitive if for any two edges $\{u,v\}, \{x,y\} \in E(X)$, there exists an automorphism $g \in \Aut(X)$ such that $\{u^g, v^g\} = \{x,y\}$.  In other words, $X$ is edge-transitive iff the action of $\Aut(X)$ on the edge set $E(X)$ has a single orbit.  

Given a group $H$ and a subset $S$ of $H$ such that $1 \notin S$ and $S = S^{-1}$, the Cayley graph of $H$ with respect to $S$, denoted by $\Cay(H,S)$, is the graph with vertex set $H$ and edge set $\{ \{h, sh\}: h \in H, s \in S \}$. The automorphism group of a Cayley graph $\Cay(H,S)$ contains the right regular representation $R(H)$ as a subgroup, whence all Cayley graphs are vertex-transitive (cf. \cite{Biggs:1993}).  Let $S$ be a set of transpositions in the symmetric group $S_n$.  The transposition graph of $S$, denoted by $T(S)$, is the graph with vertex set $[n] = \{1,\ldots,n\}$, and with vertices $i$ and $j$ being adjacent in $T(S)$ whenever $(i,j) \in S$.  A set $S$ of transpositions in $S_n$ generate $S_n$ if and only if the transposition graph $T(S)$ is connected (cf. \cite{Godsil:Royle:2001}).   

If $S$ is a set of transpositions in $S_n$, then the Cayley graph $\Cay(S_n,S)$ is called a Cayley graph generated by transpositions.  The family of Cayley graphs generated by transpositions has been well-studied because it is a suitable topology for consideration in interconnection networks (cf. \cite{Heydemann:1997}, \cite{Lakshmivarahan:etal:1993} for surveys). This family of graphs has better degree-diameter properties than the hypercube \cite{Akers:Krishnamurthy:1989}.  The automorphism group of Cayley graphs generated by transpositions has also been determined in some cases (cf. \cite{Feng:2006}, \cite{Ganesan:DM:2013}, \cite{Ganesan:JACO}, \cite{Zhang:Huang:2005}). In the present note, we further study the symmetry properties of $\Cay(S_n,S)$, especially with regards to how symmetry properties of $\Cay(S_n,S)$ depend on the properties of the generating set $S$. 

The main result of this note is the following:

\begin{Theorem} \label{thm:main:statement} Let $n \ge 5$. 

(a) Let $S, S'$ be sets of transpositions generating $S_n$. Then, the Cayley graphs $\Cay(S_n,S)$ and $\Cay(S_n,S')$ are isomorphic if and only if the transposition graphs $T(S)$ and $T(S')$ are isomorphic. 

(b) Let $S$ be a set of transpositions generating $S_n$.  Then, the Cayley graph $\Cay(S_n,S)$ is edge-transitive if and only if the transposition graph $T(S)$ is edge-transitive. 
\end{Theorem}

\begin{Remark}  Three comments and corollaries of  Theorem~\ref{thm:main:statement}: 

1. The reverse implication of Theorem~\ref{thm:main:statement}(a) is proved in \cite[Theorem 4.5]{Lakshmivarahan:etal:1993}.  Parts of Theorem~\ref{thm:main:statement} are stated in Heydemann et al \cite{Heydemann:etal:1999} and Heydemann \cite{Heydemann:1997} without a proof; they attribute the result to unpublished reports.  We could not find a proof of Theorem~\ref{thm:main:statement} in the literature.  

2. If the transposition graph $T(S)$ is the path graph on $n$ vertices, then the Cayley graph $\Cay(S_n,S)$ is called the bubble-sort graph of dimension $n$.  Some of the literature (cf. \cite{Konstantinova:2012} \cite{Latifi:Srimani:1995} \cite{Latifi:Srimani:1996} ) incorrectly assumes the bubble-sort graph is edge-transitive.  Since the path graph is not edge-transitive, Theorem~\ref{thm:main:statement}(b) implies that the bubble-sort graph is not edge-transitive. 

On the other hand, if $T(S)$ is the complete graph $K_n$, the cycle $C_n$ or the star $K_{1,n-1}$, then the corresponding Cayley graphs $\Cay(S_n,S)$, which are referred to as the complete transposition graph, the modified bubble-sort graph and the star graph, respectively, are edge-transitive because $K_n$, $C_n$ and $K_{1,n-1}$ are edge-transitive. 

3. The vertex-connectivity of a connected graph $X$, denoted by $\kappa(X)$, is the minimal number of vertices whose removal disconnects the graph (cf. \cite{Bollobas:1998}).  Clearly, $\kappa(X)$ is at most the minimum degree $\delta(X)$.  By Menger's theorem \cite{Menger:1927}, graphs with high connectivity have a large number of parallel paths between any two nodes, making communication in such interconnection networks efficient and fault-tolerant.  Latifi and Srimani \cite{Latifi:Srimani:1995} \cite{Latifi:Srimani:1996} proved that the complete transposition graphs have vertex-connectivity equal to the minimum degree.  

Watkins \cite{Watkins:1970} proved that the vertex-connectivity of a connected edge-transitive graph is maximum possible.  Thus, Theorem~\ref{thm:main:statement}(b) (in conjunction with the theorem of Watkins \cite{Watkins:1970}) gives another proof that many families of graphs, including the complete transposition graphs, modified bubble-sort graphs and the star graphs, have vertex-connectivity that is maximum possible.  

Incidentally,  Mader \cite{Mader:1970} showed that if $X$ is a connected vertex-transitive graph that does not contain a $K_4$, then $X$ has vertex-connectivity equal to its minimum degree.  Since all Cayley graphs generated by transpositions are bipartite, they do not contain a $K_4$, and so all connected Cayley graphs generated by transpositions have vertex-connectivity maximum possible.  This gives an independent proof of the optimal vertex-connectivity of connected Cayley graphs generated by transpositions. 
\qed
\end{Remark}

\section{Preliminaries}

Let $X = (V,E)$ be a graph.  The line graph of $X$, denoted by $L(X)$, is the graph with vertex set $E$, and $e, f \in E(X)$ are adjacent vertices in $L(X)$ iff $e,f$ are incident edges in $X$.  If two graphs are isomorphic, then clearly their line graphs are isomorphic.  A natural question is the following: if $X$ and $Y$ are connected graphs with isomorphic line graphs, are $X$ and $Y$ also isomorphic?  Whitney \cite{Whitney:1932} showed that the answer is in the affirmative, unless one of $X$ or $Y$ is $K_3$ and the other is $K_{1,3}$.  Every automorphism of a graph induces an automorphism of the line graph.   Whitney \cite{Whitney:1932} showed that we can go in the reverse direction: if $T$ is a connected graph on 5 or more vertices, then every automorphism of the line graph $L(T)$ is induced by a unique automorphism of $T$.   

\begin{Theorem} \label{thm:Whitney:graph:linegraph:sameautgroup} (Whitney \cite{Whitney:1932}, Sabidussi \cite{Sabidussi:1961}) 
Let $T$ be a connected graph on 5 or more vertices. Then, every automorphism of the line graph $L(T)$ is induced by a unique automorphism of $T$, and the automorphism group of $T$ and of $L(T)$ are isomorphic.
\end{Theorem}

In the sequel, we shall refer to the following result due to Feng \cite{Feng:2006} and its proof (the proof sketch is given below). 

\begin{Theorem} \label{thm:Feng:Aut:Sn:S:equals:AutTS}(Feng \cite{Feng:2006}) 
 Let $S$ be a set of transpositions in $S_n$ ($n \ge 3$).  Then, the group of automorphisms of $S_n$ that fixes $S$ setwise is isomorphic to the automorphism group of the transposition graph of $S$, i.e., $\Aut(S_n,S) \cong \Aut(T(S))$. 
\end{Theorem}

\noindent \emph{Proof sketch:}
In the proof of this result, the bijective correspondence between $\Aut(S_n,S)$ and $\Aut(T(S))$ is as follows.  If $g \in S_n$ is an automorphism of the transposition graph $T(S)$, then conjugation by $g$, denoted by $c_g$, is the corresponding element in $\Aut(S_n,S)$.  In the other direction, every element in $\Aut(S_n,S)$ is conjugation $c_g$ by some element $g \in S_n$, and it can be shown that if $c_g \in \Aut(S_n,S)$, then $g \in \Aut(T(S))$. 
\qed

We shall also refer to the following result. 

\begin{Proposition} \label{prop:Ge:restrictedtoS:is:in:AutLT} (Ganesan \cite{Ganesan:JACO})
Let $S$ be a set of transpositions generating $S_n$ ($n \ge 5$) and let $G$ be the automorphism group of $X = \Cay(S_n,S)$. Let $g \in G_e$. Then, the restriction map $g|_S$ is an automorphism of the line graph of the transposition graph of $S$.
\end{Proposition}

The proof can be found in \cite{Ganesan:JACO}. 
\section{Proof of Theorem~\ref{thm:main:statement}}

In this section, we prove both parts of Theorem~\ref{thm:main:statement}. 


\begin{Theorem}
 Let $S, S'$ be sets of transpositions generating $S_n$ ($n \ge 5$).  Then the Cayley graphs $\Cay(S_n,S)$ and $\Cay(S_n,S')$ are isomorphic if and only if the transposition graphs $T(S)$ and $T(S')$ are isomorphic.
\end{Theorem}

\noindent \emph{Proof}:
Let $X = \Cay(S_n,S)$ and $X' = \Cay(S_n,S')$.  Suppose $f$ is an isomorphism from the transposition graph $T(S)$ to the transposition graph $T(S')$.  We show that the Cayley graphs $X$ and $X'$ are isomorphic. Suppose $f$ takes $i$ to $i'$, for $i \in [n]$.  Since $f$ preserves adjacency and nonadjacency,  the transposition $(i,j) \in S$ iff $(i',j') \in S'$. Let $\sigma$ be the map from $S_n$ to itself obtained by conjugation by $f$.  Denote the image of $g \in S_n$ under the action of $\sigma$ by $g'$. Since $f$ is an isomorphism, it takes the edge set of $T(S)$ to the edge set of $T(S')$.  Hence, $S^\sigma = S$. 

We show that $\sigma: V(X) \rightarrow V(X')$ is an isomorphism from $X$ to $X'$. Suppose vertices $g,h$ are adjacent in $X$. Then there exists an $s \in S$ such that $sg=h$. Applying $\sigma$ to both sides, we get that $(sg)^\sigma = h^\sigma$, whence $s'g' = h'$.  Note that $s' \in S'$.  Hence, vertices $g'$ and $h'$ are adjacent in $X'$. By applying $\sigma^{-1}$ to both sides, we get the converse that if $g',h'$ are adjacent vertices in $X'$, then $g,h$ are adjacent vertices in $X$.   We have shown that $X$ and $X'$ are isomorphic. 

Now suppose the Cayley graphs $X$ and $X'$ are isomorphic, and let $f: V(X) \rightarrow V(X')$ be an isomorphism. Since $X'$ admits the right regular representation $R(S_n)$ as a subgroup of automorphisms, if $f$ takes the identity vertex $e \in V(X)$ to $h' \in V(X')$, then $f$ composed with $r_{h'}^{-1} \in R(S_n)$ takes $e$ to $e$.  Therefore, we may assume without loss of generality that the isomorphism $f$ maps the identity vertex of $X$ to the identity vertex of $X'$.  The neighbors of $e$ in the Cayley graphs $X$ and $X'$ are $S$ and $S'$, respectively. Hence, $f$ takes $S$ to $S'$. Consider the restriction map $f|_S$.  By the proof of Proposition~\ref{prop:Ge:restrictedtoS:is:in:AutLT}, the restriction map is an isomorphism from the line graph of $T(S)$ to the line graph of $T(S')$.  Denote these two transposition graphs $T(S), T(S')$ by $T, T'$, respectively, and their line graphs by $L(T), L(T')$, respectively.  We have just shown that the line graphs $L(T)$ and $L(T')$ are isomorphic. 

Since $S, S'$ generate $S_n$, their transposition graphs $T, T'$, respectively, are connected.  Because $X$ and $X'$ are isomorphic, $|E(T)| = |E(T')|$ and $|V(T)| = |V(T')|$.  Therefore, it is not possible that one of $T, T'$ is $K_3$ and the other $K_{1,3}$.  Since the line graphs $L(T)$ and $L(T')$ are isomorphic, by Whitney's Theorem~\ref{thm:Whitney:graph:linegraph:sameautgroup},  the transposition graphs $T$ and $T'$ are isomorphic. 
\qed

\begin{Proposition} \label{prop:Ge:Le:TS}
 Let $S$ be a set of transpositions generating $S_n$ ($n \ge 5$).  Let $G$ be the automorphism group of $X = \Cay(S_n,S)$ and let $L_e$ denote the set of element in $G_e$ that fixes the vertex $e$ and each of its neighbors.  Then, $G_e = L_e \rtimes \Aut(S_n,S)$. 
\end{Proposition}

\noindent \emph{Proof}:
Let $g \in G_e$. Then $g |_S$ is an automorphism of the line graph of $T(S)$ (cf. Proposition~\ref{prop:Ge:restrictedtoS:is:in:AutLT}). By Whitney's  Theorem~\ref{thm:Whitney:graph:linegraph:sameautgroup}, the automorphism $g|_S$ of the line graph of $T(S)$ is induced by an automorphism $h$ of $T(S)$. Conjugation by $h$, denoted by $c_h$, which is an element of $\Aut(S_n,S)$, has the same action on $S$ as $g$, i.e., $g|_S = c_h|_S$.  This implies that $g c_h^{-1} \in L_e$, whence $g \in L_e c_h$.  It follows that $G_e$ is contained in $L_e \Aut(S_n,S)$.  Clearly $L_e \Aut(S_n,S)$ is contained in $G_e$. Hence, $G_e = L_e \Aut(S_n,S)$.  

Since $L_e$ is a normal subgroup of $G_e$ (cf. \cite{Biggs:1993}), it remains to show that $L_e \cap \Aut(S_n,S) = 1$.   Each element in $L_e$ fixes $X_1(e)$ pointwise. The only element in $\Aut(S_n,S)$ which fixes $X_1(e)$ pointwise is the trivial permutation of $S_n$ because if $g \in \Aut(S_n,S)$ fixes $X_1(e)$ pointwise, then the restriction map $g|_S$ is a trivial automorphism of the line graph of $T(S)$, and hence is induced by the trivial automorphism $h$ of $T(S)$.  Since $g$ is conjugation by $h$ (cf. proof of Theorem~\ref{thm:Feng:Aut:Sn:S:equals:AutTS}), $g=1$.  We have shown that the only element in $\Aut(S_n,S)$ which fixes $S$ pointwise is the trivial permutation of $S_n$. It follows that $L_e \cap \Aut(S_n,S)=1$. 
\qed

A graph $X=(V,E)$ is said to be arc-transitive if for any two ordered pairs $(u,v), (x,y)$ of adjacent vertices, there is an automorphism $g \in \Aut(X)$ such that $u^g=x$ and $v^g=y$.

\begin{Theorem} \label{thm:edge:trans}
 Let $S$ be a set of transpositions generating $S_n$ ($n \ge 5$). Then, the Cayley graph $\Cay(S_n,S)$ is edge-transitive if and only if the transposition graph $T(S)$ is edge-transitive.
\end{Theorem}

\noindent \emph{Proof}:
Suppose the transposition graph $T(S)$ is edge-transitive.  Let $G$ be the automorphism group of $X = \Cay(S_n,S)$.  To prove $X$ is edge-transitive, it suffices to show that $G_e$ acts transitively on $X_1(e)$.  Let $t, k \in X_1(e) = S$.  Note that $t,k$  are edges of $T(S)$. By hypothesis, there exists an an automorphism $g \in S_n$ of $T(S)$ that takes edge $t$ to edge $k$.  Conjugation by $g$, denoted by $c_g$, is an automorphism of $S_n$ that takes permutation $t \in S_n$ to $k$.  Also, $c_g \in \Aut(S_n,S)$ (cf. proof of Theorem~\ref{thm:Feng:Aut:Sn:S:equals:AutTS}). 
Since $\Aut(S_n,S) \le G_e$, $G_e$ contains an element $c_g$ which takes $t$ to $k$. It follows that $G_e$ acts transitively on $X_1(e)$.

For the converse, suppose the Cayley graph $\Cay(S_n,S)$ is edge-transitive.  Fix $t \in S$. Let $r_t$ be the map from $S_n$ to itself that takes $x$ to $xt$.  Observe that $r_t$ takes the arc $(e,t)$ to the arc $(t,e)$ since $t^2=e$.  Hence, the Cayley graph $\Cay(S_n,S)$ is arc-transitive.  This implies that $G_e$ acts transitively on $X_1(e) = S$.   Since $G_e$ acts transitively on $X_1(e)$ and $L_e$ fixes $X_1(e)$ pointwise, the formula $G_e = L_e \Aut(S_n,S)$ (cf. Proposition~\ref{prop:Ge:Le:TS}) implies that $\Aut(S_n,S)$ acts transitively on $X_1(e)$.  

Let $t, k$ be two edges of the transposition graph $T(S)$, so that $t, k \in X_1(e)$.  By the argument in the previous paragraph, there exists an element $g \in \Aut(S_n,S)$ that takes vertex $t$ of $X$ to vertex $k$ of $X$.  By the bijective correspondence between $\Aut(S_n,S)$ and $\Aut(T(S))$ (cf. proof of Theorem~\ref{thm:Feng:Aut:Sn:S:equals:AutTS}), there exists an automorphism $h$ of $T(S)$ such that $g=c_h$, where $c_h$ denotes conjugation by $h$, and such that $h$ takes edge $t$ of the transposition graph to edge $k$. Thus, the set of permutations $\{h \in S_n: c_h \in \Aut(S_n,S)\}$ is contained in $\Aut(T(S))$ and  acts transitively on the edges of $T(S)$.
\qed
%

{
\bibliographystyle{plain}
\bibliography{refsaut}
}
\end{document}